\documentclass[twoside,11pt]{article}
\title{} \author{} \date{}
\usepackage{latexsym,amssymb,times}
\input amssym.def
\newtheorem{te}{Theorem}[section]
\newtheorem{prop}[te]{Proposition}

\newtheorem{fac}[te]{Fact}
\newtheorem{cla}[te]{Claim}

\newtheorem{rem}[te]{Remark}
\newtheorem{ex}[te]{Example}

\def\dok{\noindent{\bf Proof. }}
\def\kdok{\hfill $\Box$ \par \vspace*{2mm} }

\let\strokel\l

\def\f{\varphi}
\def\p{\psi}
\def\o{\omega}
\def\k{\kappa}
\def\l{\lambda}

\def\r{\rho}
\def\s{\sigma}

\def\z{\zeta}

\def\F{{\mathbb F}}

\def\P{{\mathbb P}}
\def\Q{{\mathbb Q}}

\def\N{{\mathbb N}}
\def\X{{\mathbb X}}
\def\Y{{\mathbb Y}}
\def\Z{{\mathbb Z}}

\def\A{{\mathbb A}}

\def\BL{{\mathbb L}}

\def\BK{{\mathbb K}}
\def\BT{{\mathbb T}}

\def\I{{\mathcal I}}
\def\O{{\mathcal O}}

\def\la{\langle}
\def\ra{\rangle}

\def\Nrev{(\N ^\N)_{\mathrm{rev}}}

\def\dom{\mathop{\mathrm{dom}}\nolimits}

\def\id{\mathop{\mathrm{id}}\nolimits}

\def\Card{\mathop{\rm Card}\nolimits}

\def\Iso{\mathop{\rm Iso}\nolimits}
\def\Aut{\mathop{\rm Aut}\nolimits}

\def\Sur{\mathop{\rm Sur}\nolimits}
\def\Sym{\mathop{\rm Sym}\nolimits}

\def\ar{\mathop{\rm ar}\nolimits}

\def\Int{\mathop{\rm Int}\nolimits}

\def\Mono{\mathop{\rm Mono}\nolimits}

\def\Cond{\mathop{\rm Cond}\nolimits}

\def\rfm{\mathop{\rm RFM}\nolimits}
\def\rc{\mathop{\rm RC}\nolimits}
\def\ru{\mathop{\rm RU}\nolimits}

\begin{document}

\thispagestyle{plain}
\begin{center}
           {\large \bf \uppercase{Reversible sequences of cardinals, reversible\\[1mm]
                                  equivalence relations, and similar structures}}%
\end{center}
\begin{center}
{\bf Milo\v s S.\ Kurili\'c\footnote{Department of Mathematics and Informatics, Faculty of Science, University of Novi Sad,
              Trg Dositeja Obradovi\'ca 4, 21000 Novi Sad, Serbia.
              email: milos@dmi.uns.ac.rs}
and Nenad Mora\v ca\footnote{Department of Mathematics and Informatics, Faculty of Science, University of Novi Sad,
              Trg Dositeja Obradovi\'ca 4, 21000 Novi Sad, Serbia.
              email:
              nenad.moraca@dmi.uns.ac.rs}}
\end{center}
\begin{abstract}
\noindent
A relational structure $\X$ is said to be reversible iff every bijective endomorphism $f:X\rightarrow X$ is an automorphism.
We define a sequence of non-zero cardinals $\la \k _i :i\in I\ra$ to be reversible
iff each surjection $f :I\rightarrow I$ such that
$\k _j =\sum _{i\in f^{-1}[\{ j \}]}\k_i$, for all $j\in I $, is a bijection, and characterize such sequences:
either $\la \k _i :i\in I\ra$ is a finite-to-one sequence, or $\k _i\in \N$, for all $i\in I$,
$K:=\{ m\in \N : \k _i =m ,\mbox{ for infinitely many } i\in I \}$
is a non-empty independent set,
and $\gcd (K)$ divides at most finitely many elements of the set $\{ \k _i :i\in I \}$.
We
isolate a class of binary structures such that a structure from the class is reversible iff
the sequence of cardinalities of its connectivity components is reversible.
In particular, we characterize reversible equivalence relations, reversible posets which are
disjoint unions of cardinals $\leq \o$, and some similar structures. In addition, we show that
a poset with linearly ordered connectivity components is reversible, if the corresponding
sequence of cardinalities is reversible and, using this fact, detect a wide class of examples of reversible posets and topological spaces.

{\sl 2010 MSC}:
03C50, 
03C07, 
03E05, 
06A06, 
05C20, 
05C40. 

{\sl Key words}: reversible sequence of cardinals, reversible sequence of natural numbers, reversible equivalence relation, digraph, poset.
\end{abstract}
\section{Introduction}\label{S1}
A structure is called reversible iff all its bijective endomorphisms are automorphisms
and the class of reversible structures contains, for example, Euclidean, compact and many other relevant topological spaces \cite{RajWil,DoyHoc,Dow},
linear orders, Boolean lattices, well founded posets with finite levels \cite{Kuk,Kuk1},
tournaments, Henson graphs \cite{KuMoExtr}, and Henson digraphs \cite{KRet}.
In addition, reversible structures have several distinguished properties; for example, the Cantor-Schr\"{o}der-Bernstein property for condensations
(bijective homomorphisms).

It seems that the property of reversibility of relational structures is more of set-theoretical or combinatorial,
than of model-theoretical nature--it is an invariant of isomorphism and condensational equivalence, while it is not preserved under
bi-embeddability, bi-definability and elementary equivalence \cite{KDef,KuMoSim}.
But it is an invariant of some forms of bi-interpretability \cite{KRet},
extreme elements of $L_{\infty \o}$-definable classes of structures are reversible under some syntactical restrictions \cite{KuMoExtr},
and all structures first-order definable in linear orders by quantifier-free formulas without parameters (i.e., monomorphic or chainable structures)
are reversible \cite{KDef}.

In this article we continue the investigation of reversibility in the class of disconnected binary structures initiated in \cite{KuMoDiscI}.
If $\X$ is a binary structure and $\X _i$, $i\in I$, are its connectivity components, then, clearly, the sequence of cardinal numbers $\la |X_i|:i\in I\ra$
is an isomorphism-invariant of the structure and in some classes of structures (for example, in the class of equivalence relations)
that cardinal invariant characterizes  the structure up to isomorphism.
In such classes the reversibility of a structure, being an isomorphism-invariant as well,
can be regarded as a property of the corresponding sequence of cardinals.

So, using the characterization of reversible disconnected binary structures from \cite{KuMoDiscI} (see Fact \ref{TB044})
we easily isolate the following property of sequences of cardinals (called reversibility as well)
which characterizes reversibility in the class of equivalence relations:
If $I$ is a non-empty set, an $I$-sequence of non-zero cardinals $\la \k _i :i\in I\ra$ will be called
{\it reversible} iff there is no non-injective surjection $f :I\rightarrow I$ such that
\begin{equation}\label{EQB032}\textstyle
\forall j\in I \;\;\k _j =\sum _{i\in f^{-1}[\{ j \}]}\k_i .
\end{equation}
The first main result of this paper is the following characterization of reversible sequences of cardinals.
In order to state it we recall some definitions.
For a  subset $K$ of the set of natural numbers, $\N$, let $\la K\ra$ denote the subsemigroup of the semigroup  $\la \N , +\ra$
generated by $K$. A set $K$ is called {\it independent} iff
\begin{equation}\label{EQB003}
\forall n\in K \;\; n\not\in \la K\setminus \{ n \}\ra.
\end{equation}
So, $\emptyset$ is an independent set. If $K \neq\emptyset$, by $\gcd (K)$ we denote the greatest common divisor of the numbers from $K$.
\begin{te}\label{TB042}
A sequence of non-zero cardinals $\la \k _i :i\in I\ra$ is reversible iff
\begin{itemize}
\item[-] either $\la \k _i :i\in I\ra$ is a finite-to-one sequence,
\item[-] or $\k _i\in \N$, for all $i\in I$,\\
$K:=\{ m\in \N : |\{ i\in I : \k _i =m\}|\geq \o\}$
is a non-empty independent set, \\
and $\gcd (K)$ divides at most finitely many elements of the set $\{ \k _i :i\in I \}$.\footnote{For example, if $I$ is a non-empty set of any size
and $\la n_i :i\in I\ra \in {}^{I}\N$, then by Theorem \ref{TB042} we have:
 if $K=\emptyset$ (which is possible if $|I|\leq \o$), then $\la n_i \ra $ is a reversible sequence;
 if $K=\{ 2,5\}$, then  $\la n_i \ra $ is a reversible sequence iff the set $\{ n_i :i\in I \}$ is finite;
 if $K=\{ 4,10\}$, then  $\la n_i \ra $ is a reversible sequence iff the set $\{ n_i :i\in I \}$
contains at most finitely many even numbers.
}
\end{itemize}
\end{te}
A proof of Theorem \ref{TB042} is given in the last (and the largest) Section \ref{S4}, where,
in addition, we show that the set of reversible sequences of natural numbers is a dense $F_{\s \delta \s}$-subset of the Baire space,
and that it is not a subsemigroup of $\la \N ^\N ,\circ\ra$.

Section \ref{S2} contains definitions and facts making the paper  self-contained.

In Section \ref{S3}, generalizing the situation with equivalence relations, we isolate a wider class of
structures with the same property--that the reversibility of a structure from the class is equivalent
to the reversibility of the corresponding sequence of sizes of its components--the class of structures having
the sequence of components rich for monomorphisms. We also study the class $\rfm$ of such sequences of structures,
compare it with some relevant classes, detect some classes of structures such that the reversibility of a structure from the class
follows from the reversibility of the corresponding cardinal sequence and in this way detect
wide classes of reversible digraphs, posets, and topological spaces.
\section{Preliminaries}\label{S2}
\paragraph{Reversible structures}
If $L=\la R_i :i\in I\ra$ is a relational language, where $\ar (R _i)=n_i\in \N$, for $i\in I$, and  $\X$ and $\Y$ are $L$-structures, then by $\Iso (\X ,\Y )$, $\Cond (\X ,\Y )$ and $\Mono (\X ,\Y )$
we denote the set of all isomorphisms, condensations (bijective homomorphisms) and monomorphisms (injective homomorphisms) from $\X$ to $\Y$ respectively.
Clearly, $\Iso (\X ,\X )$ is the set of automorphisms, $\Aut (\X )$, of $\X$, instead of $\Cond (\X ,\X )$ we will write $\Cond (\X )$ etc.
For a set $X$ by $\Sym (X)$ (resp.\ $\Sur (X)$) we denote the set of all bijections (resp.\ surjections) $f:X\rightarrow X$.

The {\it condensational preorder} $\preccurlyeq _c $ on the class of $L$-structures is defined by $\X \preccurlyeq _c \Y$ iff $\Cond (\X ,\Y )\neq\emptyset$,
the {\it condensational equivalence} is the equivalence relation defined on the same class by $\X \sim _c \Y$ iff
$\X \preccurlyeq _c \Y$ and $\Y \preccurlyeq _c \X$
and it determines the antisymmetric quotient of the condensational preorder, the {\it condensational order}, in the usual way.

An $L$-structure $\X =\la X,\r\ra$ is called {\it reversible} iff $\Cond (\X )=\Aut (\X )$.
Clearly, $\r=\la \r _i :i\in I\ra$ is an element of the set $\Int _L (X)=\prod _{i\in I}P(X^{n_i})$
of all interpretations of the language $L$ over the domain $X$ and defining the partial order $\subset$ on $\Int _L (X)$ by
$\r\subset \s$ iff $\r_i \subset \s _i$, for all $i\in I$, it is easy to obtain the following simple characterizations
of reversible $L$-structures (see \cite{KuMoVar}).
\begin{fac}\label{TB045}
For an $L$-structure $\X =\la X,\r\ra$ the following conditions are equivalent

(a) $\X $ is a reversible structure,

(b) $\forall\s\in\Int _L (X)\;\;(\s\varsubsetneq \r \Rightarrow \s\not\cong \r )$,

(c) $\forall\s\in\Int _L (X)\;\;(\r\varsubsetneq \s \Rightarrow \s\not\cong \r )$,

(d) $\forall f\in \Sym (X)\;\; ( f[\r ]\subset \r \Rightarrow f[\r ]= \r)$.
\end{fac}
Reversible $L$-structures
have the Cantor-Schr\"{o}der-Bernstein property for condensations. Moreover we have (see \cite{KuMoVar})
\begin{fac}\label{TA011}
Let $\X$ and $\Y$ be $L$-structures. If $\X$ is a reversible structure and $\Y \sim _c \X$, then
$\Y \cong \X$ (thus $\Y$ is reversible too) and $\Cond (\X ,\Y )=\Iso (\X ,\Y )$.
\end{fac}
\paragraph{Disconnected binary structures}
Let $L_b$ be the binary language, that is,
$L_{b}=\la R\ra$ and $\ar (R)=2$. If $\X=\la X,\rho  \ra$ is an $L_b$-structure, then
the transitive closure $\rho  _{rst}$ of the relation $\rho  _{rs} =\Delta _X \cup \rho  \cup \rho  ^{-1}$ (given by $x \,\rho _{rst} \,y$ iff there are $n\in \N$ and
$z_0 =x , z_1, \dots ,z_n =y$ such that $z_i \;\rho  _{rs} \;z_{i+1}$, for each $i<n$)
is the minimal equivalence relation on $X$ containing $\rho $. The corresponding equivalence classes are called the {\it components} of $\X$ and
the structure $\X$ is called {\it connected} iff $|X/\rho  _{rst} |=1$.

If $\X _i=\la X_i, \rho  _i \ra$, $i\in I$, are connected $L_b$-structures  and $X_i \cap X_j =\emptyset$, for
different $i,j\in I$, then the structure $\bigcup _{i\in I} \X _i =\la \bigcup _{i\in I} X_i , \bigcup _{i\in I} \rho  _i\ra$ is
the {\it disjoint union} of the structures $\X _i$, $i\in I$,  and the structures $\X _i$, $i\in I$, are its components.
\begin{fac}[\cite{KuMoDiscI}]\label{TB044}
Let $\X _i$, $i\in I$, be pairwise disjoint and connected $L_b$-structures.
\begin{itemize}
\item[(a)] If  $\;\bigcup _{i\in I} \X _i$ is reversible, then all structures $\mathbb{X}_i$, $i\in I$, are reversible.
\item[(b)] $\bigcup _{i\in I} \X _i$ is a reversible structure iff

whenever  $f:I\rightarrow I$ is a surjection,
$g_i \in \Mono(\X _i ,\X _{f(i)})$, for $i\in I$, and
\begin{equation}\label{EQB035}
\forall j\in I \;\; \Big(\Big\{g_i[X_i]: i\in f^{-1}[\{ j\}] \Big\} \mbox{ is a partition of }X_j\Big),
\end{equation}
we have
\begin{equation}\label{EQB034}
f\in \Sym (I) \;\land \;\forall i\in I \;\; g_i \in \Iso (\X _i,\X _{f(i)}) .
\end{equation}
\end{itemize}
\end{fac}
\section{Sequences of structures rich for monomorphisms}\label{S3}
We will say that a sequence of $L$-structures $\la \X_i :i\in I\ra$ is
{\it rich for monomorphisms} iff
\begin{equation}\label{EQA027}
\forall i,j\in I \;\; \forall A\in [X_j ]^{|X_i|} \;\;\exists g\in \Mono (\X _i , \X _j)\;\; g[X_i]=A.
\end{equation}
By Fact \ref{TB044}(a), a necessary condition for the reversibility of a disconnected binary structure is the reversibility of its components.
Hence, and in order to simplify notation, in the sequel we work under the following assumption:
\begin{itemize}
\item[($\ast$)] $\X _i$, $i\in I$, are pairwise disjoint, connected and reversible $L_b$-structures.
\end{itemize}
Let $\rfm$ denote the class of sequences of $L_b$-structures $\la \X_i :i\in I\ra$ (where $I$ is any non-empty set)
satisfying $(\ast)$ and which are rich for monomorphisms.
\subsection{Reversible equivalence relations and similar structures}
First we show that the reversibility of a structure having the sequence of components in $\rfm$ depends only on the corresponding cardinal sequence.
\begin{te}\label{TA021}
If $\la \X_i :i\in I\ra \in \rfm$, then

(a) The structures of the same size are isomorphic, 

(b) $\bigcup _{i\in I}\X_i$ is reversible $\Leftrightarrow$  $\la |X_i|: i\in I \ra$ is a reversible sequence of cardinals.
\end{te}
\dok
(a) If $|X_i|=|X_j|$, then by (\ref{EQA027}) there are $g\in \Cond (\X _i , \X _j)$ and $g'\in \Cond (\X _j , \X _i)$. So
$\X _i \sim _c \X _j$, which by Fact \ref{TA011} implies that $\X _i \cong \X _j$.

(b) ($\Rightarrow$) Suppose that the sequence $\la |X_i|: i\in I \ra$ is not reversible and that
$f:I\rightarrow I$ is a noninjective surjection such that for each $j\in I$ we have $|X_j| =\sum _{i\in f^{-1}[\{ j \}]}|X_i|$.
Then for $j\in I$ there is a partition $\{ A^j_i:i\in f^{-1}[\{ j \}]\}$ of $X_j$ such that $|A^j_i|=|X_i|$, for all $i\in f^{-1}[\{ j \}]$
and, by (\ref{EQA027}), there are monomorphisms $g_i : \X _i \rightarrow \X_j =\X _{f(i)}$ satisfying $g_i [X_i]=A^j_i$.
By Fact \ref{TB044}(b) the structure $\bigcup _{i\in I}\X_i $ is not reversible.

($\Leftarrow$) Let $\la |X_i|: i\in I \ra$ be a reversible sequence of cardinals.
In order to use Fact \ref{TB044}(b), assuming that $f:I\rightarrow I$ is a surjection,
$g_i \in \Mono(\X _i ,\X _{f(i)})$, for $i\in I$, and that (\ref{EQB035}) holds,
we prove (\ref{EQB034}).
First, for $i\in I$, since the function $g_i$ is injection we have $|X_i|=|g_i[X_i]|$.
So, by (\ref{EQB035}) for each $j\in I$ we have $|X_j| =\sum _{i\in f^{-1}[\{ j \}]}|X_i|$ and,
since the sequence $\la |X_i|: i\in I \ra$ is reversible, $f\in \Sym (I)$.

Consequently, for $i\in I$ we have $g_i[X_i]=X_{f(i)}$ and, hence, $|X_i|=|X_{f(i)}|$,
which by (a) implies $\X _i \cong \X _{f(i)}$ and, in addition, $g_i\in \Cond (\X _i, \X _{f(i)})$.
Since the structures $\X_i$ are reversible, by Fact \ref{TA011} we have
$\Cond (\X _i, \X _{f(i)})=\Iso (\X _i, \X _{f(i)})$; so $g_i\in \Iso (\X _i, \X _{f(i)})$, for all $i\in I$,
and  (\ref{EQB034}) is true indeed.
\hfill $\Box$
\begin{te}\label{TB048}
Let  $\sim$ be an equivalence relation on a set $X$, $\X =\la X ,\sim\ra$,  and $\{ X_i :i\in I\}$ the corresponding partition.
Then the structure $\X$ is reversible iff $\la |X_i| :i\in I \ra$ is a reversible sequence of cardinals.

The same holds for the graphs (resp.\ posets) of the form $\X =\bigcup _{i\in I}\X_i$, where $\X_i$, $i\in I$, are pairwise disjoint
complete graphs (resp.\ ordinals $\leq \o$).
\end{te}
\dok
It is clear that any sequence of disjoint $L_b$-structures with full relations, or complete graphs, or well orders $\leq \o$ belongs to
$\rfm$; so Theorem \ref{TA021} applies.
\hfill $\Box$
\begin{rem}\label{RB001}\rm
There are ${\mathfrak c}$-many non-isomorphic countable reversible equivalence relations
(and the same holds for the classes of graphs and posets from Theorem \ref{TB048}).
By Theorems \ref{TB048} and \ref{TB042}, if  $\la n_i :i\in \N\ra \in {}^{\N}\N$ is an increasing sequence,
then the structure $\X_{\la n_i\ra}$ with the equivalence relation on $\N$
determined by a partition $\{ C_i :i\in \N \}$, where $|C_i|=n_i$, for all $i\in \N$, is reversible.
Also, if $\la n_i :i\in \N\ra\neq \la n'_i :i\in \N\ra $, then the corresponding structures are non-isomorphic.
For $A\in [\N ]^\o$ let $\la n^A_i :i\in \N\ra$ be the increasing enumeration of the set $A$. Then the structures
$\X _{\la n^A_i\ra}$, $A\in [\N ]^\o$, are non-isomorphic, countable and reversible.
\end{rem}
\subsection{More reversible digraphs, posets, and topological spaces}
In the following theorem we detect a class of structures such that the reversibility of a structure belonging to the class {\it follows} from the
reversibility of the sequence of cardinalities of its components.
\begin{te}\label{TB047}
If $\X_i$, $i\in I$, are disjoint tournaments and the sequence of cardinals $\la |X_i|: i\in I\ra$ is reversible,
then the digraph $\bigcup _{i\in I}\X_i$ is reversible.

This statement holds if, in particular, $\X_i$, $i\in I$, are disjoint linear orders. Then $\bigcup _{i\in I}\X_i$ is a reversible
disconnected partial order.
\end{te}
\dok
In order to apply Fact \ref{TB044}(b) we suppose that $f:I\rightarrow I$ is a surjection,
$g_i \in \Mono(\X _i ,\X _{f(i)})$, for $i\in I$, and that (\ref{EQB035}) holds.
Then, since $|g_i[X_i]|=|X_i|$, for $i\in I$, for each $j\in I$ by (\ref{EQB035}) we have $|X_j|=\sum _{i\in f^{-1}[\{ j\} ]}|X_i|$, which,
since the sequence $\la |X_i|: i\in I\ra$ is reversible, implies that $f\in \Sym (I)$.
Thus for each $i\in I$ we have $g_i \in \Cond(\X _i ,\X _{f(i)})$, and, since the structures $\X _i$, $i\in I$,
are tournaments, $\Cond(\X _i ,\X _{f(i)})= \Iso(\X _i ,\X _{f(i)})$. Thus (\ref{EQB034}) is true and the digraph $\bigcup _{i\in I}\X_i$ is reversible
indeed.
\hfill $\Box$
\begin{ex}\label{EXB011}\rm
The converse of Theorem \ref{TB047} is not true. Let $I=\N$ and
$\X _i \cong \o i$, for $i\in \N$. By Theorem \ref{TB042} the sequence of cardinals $\la \o, \o , \dots\ra$
is not reversible. Using Fact \ref{TB044}(b) we show that $\X =\bigcup _{i\in \N}\X_i$ is a reversible structure.
Let $f:\N\rightarrow \N$ be a surjection,
$g_i \in \Mono(\X _i ,\X _{f(i)})$, for $i\in \N$, and let (\ref{EQB035}) hold.
First, by induction we show that $f(i)=i$, for all $i\in \N$.

If $i\in \N$ and $f(i)=1$, then $g_i \in \Mono(\X _i , \X _1 )$ and,
since monomorphisms between linear orders are embeddings, $\o i \hookrightarrow \o$
and, hence, $i=1$. Thus $f^{-1}[\{ 1\}]\subset \{ 1\}$ and, since $f$ is a surjection, $f^{-1}[\{ 1\}]= \{ 1\}$.

Let $j\in \N$ and $f(k)=k$, for all $k<j$. If $i\in \N$ and $f(i)=j$, then $g_i \in \Mono(\X _i , \X _j )$ and,
as above, $\o i \hookrightarrow \o j$, which means that $i\leq j$. By the induction hypothesis we have $i\geq j$,
so $i=j$ and, thus, $f^{-1}[\{ j\}]\subset \{ j\}$ and, since $f$ is a surjection, $f^{-1}[\{ j\}]= \{ j\}$.

So, $f=\id _\N\in \Sym (\N )$, which by (\ref{EQB035}) implies that for each $i\in \N$ we have $g_i \in \Cond(\X _i ,\X _i)=\Iso(\X _i ,\X _i)$
and (\ref{EQB034})  is proved.
\end{ex}
\begin{ex}\label{EXB015}\rm
More reversible posets and topological spaces.
The reversible posets constructed in Examples \ref{EXB011} and \ref{EXB014} are well-founded and with infinite levels.
More generally, by Theorem \ref{TB047}, if $\la \k _i : i\in I\ra$ is {\it any} reversible sequence of cardinals
(e.g., if it is finite-to-one, if we would like infinite components) and $L_i$, $i\in I$, are {\it any} linear orders, where $|L_i|=\k _i$, then the poset $\bigcup _{i\in I}L_i$  is reversible.

Recalling that if $\P =\la P, \leq \ra$ is a partial order and $\O$ the topology on the set $P$
generated by the base consisting of the sets of the form $B_p:= \{ q\in p: q\leq p\}$, then endomorphisms of $\P$
are exactly the continuous self mappings of the space $\la P,\O\ra$, we conclude that
the poset $\P$ is reversible iff $\la P,\O\ra$ is a reversible topological space (i.e., each continuous bijection is an automorphism).
So, Examples \ref{EXB011},  \ref{EXB014} and Theorem \ref{TB047} generate a large class of reversible topological spaces.
\end{ex}
\subsection{More sequences from RFM}
We recall that a relational structure $\X$ is called {\it monomorphic}
iff each two finite substructures of $\X$ of the same size are isomorphic,
and that, by the well-known theorems of Fra\"{\i}ss\'{e} (for finite languages) and Pouzet (for languages and structures of any size), see \cite{Fra},
an infinite structure $\X$ is monomorphic iff it is {\it chainable} i.e.\ there is a linear order $\prec$ on its domain, $X$,
such that the relations of $\X$ are definable in the structure $\la X,\prec\ra$ by quantifier-free formulas without parameters.
Then it is said that $\prec$ {\it chains} $\X$, or that $\X$ is {\it chainable} by $\prec$.
For convenience, a structure $\X$ will be called {\it copy-maximal} (resp.\ {\it mono-range-maximal}) iff
for each $A\in [X ]^{|X|}$ there is an embedding (resp.\ a monomorphism) $g:\X \rightarrow \X$ satisfying $g[X]=A$.

By (\ref{EQA027}),  Theorem \ref{TA021}(a) and since each set of cardinals is well ordered, 
a sequence $\la \X _i :i\in I\ra\in \rfm$ can be described in the following way.
There are an ordinal $\eta$ and a sequence of connected reversible $L_b$-structures $\la \Y _\xi :\xi <\eta\ra$ (the {\it range}) such that, defining $\k _\xi := |Y_\xi|$, we have

(r1) $\xi < \zeta <\eta \Rightarrow \k_\xi<\k _\zeta$,

(r2) $\Y _\xi$ is a mono-range-maximal structure, for each $\xi <\eta$,

(r3) $\xi < \zeta <\eta \Rightarrow \forall A\in [Y_\zeta ]^{\k _\xi } \;\Cond (\Y _\xi ,A)\neq\emptyset$,

\noindent
and there is a surjection $h:I\rightarrow \eta$ such that for each $\xi <\eta $ and $i\in h^{-1}[\{ \xi \}]$ we have $\X_i\cong \Y _\xi$,
and $X_i \cap X_j =\emptyset$, for $i\neq j$.
So, by  Theorem \ref{TA021}(b), the structure
$\bigcup _{i\in I}\X _i$ is reversible iff $\la \k _{h(i)}:i\in I \ra$ is a reversible sequence of cardinals.
Here we consider conditions (r2) and (r3).
\paragraph{Condition (r2)}
Clearly, condition (r2) will be satisfied if the structures $\Y _\xi$ are finite or copy-maximal.
From more general results of Gibson, Pouzet and Woodrow \cite{Gib} it follows that a structure $\X$ of size $\k\geq \o$
is copy-maximal iff it is $\k$-chainable, that is, there is a linear order $\prec$ on $X$ which chains $\X$ and $\la X ,\prec\ra\cong \la \k ,<\ra$.
On the other hand, a simple application of  Ramsey's theorem shows that, up to isomorphism,
there are only eight countable binary copy-maximal structures and the same
holds for uncountable binary structures (see also \cite{Ktow,KZb}). The six connected of them are
$\la \k , \k ^2\ra$, $\la \k , \k ^2 \setminus \Delta _\k \ra$, $\la \k , <\ra$, $\la \k , \leq \ra$ $\la \k , >\ra$, and $\la \k , \geq \ra$, and
they are reversible.
In addition, since in the class of linear orders monomorphisms are embeddings, mono-range-maximal linear orders are copy-maximal
thus the only four mono-range-maximal linear orders of size $\k$ are mentioned above.
The following example shows that the class of mono-range-maximal posets is not so restrictive.
\begin{ex}\label{EXB009}\rm
The posets of the form $\X_{\l ,\k}:=\A_\l +\BL_\k$, where $2\leq \l <\k\geq \o$, $\A_\l$ is an antichain of size $\l$,
and $\BL _\k \cong \la \k , <\ra$, are not copy-maximal and, moreover, if $\l \geq \o$, $\X_{\l ,\k}$ is not almost chainable (see \cite{Fra, Gib} for details).
But $\X_{\l ,\k}$ is mono-range-maximal (if $S \in [X]^\k$, then $S\cong \A_\mu +\BL_\k$, for some $\mu \leq \l$, and it is easy
to construct a monomorphism from $\X_{\l ,\k}$ onto $S$). If $\l <\o$, then  $\X_{\l ,\k}$ is a well-founded poset with finite levels so, by
\cite{Kuk}, it is reversible.
\end{ex}
\paragraph{Condition (r3)}
All the structures considered in Theorem \ref{TB048} -
disjoint unions of
(a) structures with full relations, 
(b) complete graphs, and
(c) ordinals $\leq \o$,
give examples of sequences satisfying (r3) and all of them have monomorphic components. The following examples show that this condition
is not necessary for application of  Theorem  \ref{TA021}(b).
\begin{ex}\label{EXB013}\rm
Structures from $\rfm$ with non-monomorphic components.
Let

- $\BT _3$ be the three-element tree $\la \{ 0,1,2\}, \{ \la 0,1\ra , \la 0,2\ra\}\ra$,

- $\BL _5$ the five-element linear order,

- $\BK _6^*$ a complete graph with 6 nodes and 3 of them reflexified (loops),

- $\F _8$ the eight-element structure with the full relation.

\noindent
Now, if $\k$ and $\l$ are infinite cardinals, $m,n\in \o$ and $\X$ is the (pairwise disjoint) union of
$\k$-many copies of $\BT _3$, $\l$-many copies of $\BL _5$, $m$ copies of $\BK _6^*$ and $n$ copies of $\F _8$, then
the sequence $\la \BT _3  ,\BL _5,\BK _6^* ,\F _8\ra $ satisfies (r1)-(r3),
the corresponding sequence of components of $\X$ belongs to $\rfm$
and $\X$ is reversible because,
in notation of Proposition \ref{TB037}, $K=\{ 3,5\}$ and the set $\{ n_i : i\in I\}=\{3,5,6,8\}$ is finite and we apply  Theorem  \ref{TA021}(b).
\end{ex}
\begin{ex}\label{EXB014}\rm
A structure from $\rfm$ having all components non-monomorphic.
Let $\X_{2,\k}=\A _2 +\BL _\k$, for $1\leq \k \leq \o$, be the posets defined as in Example \ref{EXB009}.
It is easy to see that $\la \X _{2,\k} : 1\leq \k \leq \o\ra \in \rfm$ .
Since the corresponding sequence of cardinals
$\la 3,4,5, \dots ,\o\ra$ is one-to-one and, thus, reversible, the structure $\X =\bigcup _{1\leq \k \leq \o}\X _{2,\k}$ is reversible.
Clearly, its components, $\X _{2,\k}$, are not 2-monomorphic.
\end{ex}
\subsection{The classes RFM, RC, and RU}
If by $\rc$ (resp.\ $\ru$)
we denote the class of sequences  $\la \X_i :i\in I\ra$ 
satisfying $(\ast)$ and such that $\la |\X_i | :i\in I\ra$ is a reversible sequence of cardinals,
(resp.\  the structure $\bigcup _{i\in I}\X_i$ is reversible),
then by  Theorem \ref{TA021}(b) we have $\rfm \cap \ru =\rfm \cap \rc$. The following
example shows that this equality is the only constraint, regarding the relationship between the classes $\rfm$, $\rc$ and $\ru$.
\begin{ex}\label{EXB010}\rm
(a) $\rfm \setminus (\ru \cup \rc )\neq\emptyset$.
If $\X _i \cong \la \o , <\ra$, for $i\in \o$, then by Theorem \ref{TB042} the sequence of cardinals $\la \o, \o , \dots\ra$
is not reversible but, since ($\la A, < \,\upharpoonright \!A\ra \cong \la \o ,< \ra$, for each $A\in [\o ]^\o$,
the sequence $\la \X_i :i\in I\ra$ is rich for monomorphisms. It is easy to see that the structure $\bigcup _{i\in I}\X_i$ is not reversible.

(b) $\rc \setminus (\rfm \cup \ru)\neq\emptyset$. Let $\X = \la \Z ,\r \ra$, where $\r =\{ \la i,i \ra :i\geq 0 \}$. Then
$\X=\bigcup _{i\in \Z}\X_i$, where $\X_i=\la \{ i\} , \emptyset \ra$, for $i<0$, and $\X_i=\la \{ i\} , \{ \la i,i \ra \} \ra$, for $i\geq 0$.
The corresponding sequence of cardinals $\la \dots ,1, 1 , \dots\ra$ is reversible and, since $\X \cong \la \Z ,\r \setminus \{ \la 0,0\ra\}\ra$, by Fact \ref{TB045} the structure $\bigcup _{i\in \Z}\X_i$ is not reversible. Since $\X _{-1}\not\cong \X _0$, by  Theorem  \ref{TA021}(a)
the sequence of structures $\la \X_i :i\in \Z\ra$
is not rich for monomorphisms.

(c) $\ru \setminus (\rfm \cup \rc)\neq\emptyset$. Let $\X = \la \Z ,\r \ra$, where
$
\r =\{ \la i,i \ra :i< 0 \}\cup \{ \la 2i, 2i+1 \ra :i\geq 0\}.
$
Then we have
$\X=\bigcup _{i\in \Z}\X_i$, where $\X_i=\la \{ i\} , \{ \la i,i \ra \} \ra$, for $i<0$,
and $\X_i=\la \{ 2i, 2i+1  \} , \{ \la  2i, 2i+1 \ra \} \ra$, for $i\geq 0$.
Now, the corresponding sequence of cardinals $\la \dots ,1,1,2,2,  \dots\ra$ is not reversible, because the set $K=\{ 1,2\}$ is not independent ($1+1=2$).
Since $\Mono (\X _{-1},\X _0)=\emptyset$ we have $\la \X_i :i\in \Z\ra \not\in \rfm$.
But, by Fact \ref{TB045}, the structure $\bigcup _{i\in \Z}\X_i$ is reversible,
namely, if $\s \varsubsetneq \r$, then the structure $\la \Z ,\s \ra$ has an one-element component with the empty relation and, hence,
it is not isomorphic to $\X$.

(d) $(\ru \cap \rc)\setminus \rfm\neq\emptyset$. Let $I$ be the ordinal $\o +2 =\o \cup \{ \o ,\o +1\}$ and let $\X =\bigcup _{i\in \o +2}\X_i$,
where $\X_i$ are pairwise disjoint linear orders such that $\X_i \cong i+1$, for $i\in\o$, $\X _\o \cong \o$, and $\X _{\o +1}\cong \Q$.
The corresponding sequence of cardinals $\la 1,2,\dots , \o ,\o\ra$ is finite-to-one and, by Theorem \ref{TB042}, reversible. By Theorem
\ref{TB047} the union $\bigcup _{i\in I}\X_i$ is reversible too. Since $\o \not\cong \Q$ by  Theorem \ref{TA021}(a) we have $\la \X_i :i\in I\ra \not\in \rfm$.
\end{ex}
Let $\rfm _{LO}$, $\rc_{LO}$ and $\ru_{LO}$ denote the classes of sequences  of linear orders $\la \X_i :i\in I\ra$
belonging to classes $\rfm$, $\rc$ and $\ru$. Here, by Theorem \ref{TB047} we obtain one more constraint:
$\rc_{LO} \subset \ru_{LO}$, and the following example shows that, in general, there are no more constraints.
\begin{ex}\label{EXB012}\rm
$\rfm _{LO}\setminus \ru_{LO}\neq\emptyset$ is witnessed by the poset $\bigcup _\o \o$, from Example \ref{EXB010}(a).
The poset $\bigcup _{n\in \N}n \cup \o \cup \Q$ from Example \ref{EXB010}(d) belongs to the class $\rc_{LO}\setminus\rfm _{LO}$,
while the poset $\bigcup _{n\in \N}\o n$ (see Example \ref{EXB011}) belongs to the class $\ru_{LO}\setminus \rc_{LO}$.
\end{ex}
\section{Reversible cardinal sequences -- a proof of Theorem \ref{TB042}}\label{S4}
Theorem \ref{TB042} follows from Propositions \ref{TA020} and \ref{TB037} given in the sequel.

If $\la \k _i :i\in I\ra$ is a sequence of cardinals and $\k$ a cardinal, let
$$
I_\k := \{ i\in I : \k _i =\k\} .
$$
\subsection{Reduction to the case when the cardinals are finite}
\begin{prop}\label{TA020}
A sequence of non-zero cardinals $\la \k _i :i\in I\ra$ is reversible iff it is a finite-to-one sequence or a reversible sequence in $\N$.
\end{prop}
\dok
The implications ``$ \Leftarrow $" and ``$ \Rightarrow $" follow from Claims \ref{TB040} and \ref{TB041} respectively.
\begin{cla}\label{TB040}
If $\la \k _i :i\in I\ra$ is a finite-to-one sequence, it   is reversible.
\end{cla}
\dok
Let $|I_\k |<\o $, for all $\k \in \Card$. The {\it set} $\{ \k_i:i\in I \}$ is well-ordered
and, hence, there is an ordinal $\z$
and an enumeration  $\{ \k_i:i\in I \} =\{ \k _\xi : \xi <\z \}$
such that $ \xi <\xi'$ implies $ \k _\xi <\k _{\xi'}$.
Assuming that $f:I\rightarrow I$ is a surjection satisfying (\ref{EQB032}) we show that $f$ is a bijection.
First, by induction we prove that
\begin{equation}\label{EQA026}
\forall \xi <\z \;\; f[I _{\k _\xi}]=I _{\k _\xi}.
\end{equation}
If $j\in I_{\k _0}$, then, by (\ref{EQB032}),
for $i\in f^{-1}[\{ j \}]$ we have $\k_i \leq \k_j =\k _0$, which, by the minimality of $\k_0$,
implies that $\k_i=\k _0$, that is, $i\in I_{\k _0}$.
Thus $f^{-1}[\{ j \}]\subset I_{\k _0}$, for all $j\in I_{\k _0}$, and, hence,
$f^{-1}[I_{\k _0}]\subset I_{\k _0}$. Since $f$ is onto we have $I_{\k _0} =f[f^{-1}[I_{\k _0}]]\subset f[I_{\k _0}]$ thus
$|I_{\k _0}|\leq |f[I_{\k _0}]|\leq |I_{\k _0}|$ and, hence,  $|f[I_{\k _0}]|= |I_{\k _0}|$, which, since the set $I_{\k _0}  $
is finite and $I_{\k _0}\subset f[I_{\k _0}]$, implies that $f[I_{\k _0}]= I_{\k _0}$.

Assuming that $\eta <\z$ and $f[I _{\k _\xi}]=I _{\k _\xi}$, for all $\xi <\eta$, we prove  $f[I _{\k _\eta}]=I _{\k _\eta}$.
If $j\in I_{\k _\eta}$,
then, by (\ref{EQB032}), for $i\in f^{-1}[\{ j \}]$ we have
$\k_i\leq \k_j=\k _\eta$.
The inequality $\k_i<\k _\eta$ would imply that $\k_i=\k _\xi$, for some $\xi <\eta$,
and, hence, $i\in I _{\k _\xi}$
and, by the induction hypothesis, $f(i)=j\in I _{\k _\xi}$, which is not true.
Thus $\k_i=\k _\eta$
and, hence, $i\in I_{\k _\eta}$.
Thus $f^{-1}[\{ j \}]\subset I_{\k _\eta}$, for all $j\in I_{\k _\eta}$, and, hence,
$f^{-1}[I_{\k _\eta}]\subset I_{\k _\eta}$.
Now, as above we show that $f[I_{\k _\eta}]= I_{\k _\eta}$ and (\ref{EQA026}) is proved.

By (\ref{EQA026}) and since the sets $I_{\k _\xi}$  are finite,
the restrictions $f\upharpoonright I_{\k _\xi}: I_{\k _\xi}\rightarrow I_{\k _\xi}$, $\xi <\z$, are bijections
and, since $\{ I _{\k _\xi}: \xi <\z \}$ is a partition of the set $I$, $f$ is a bijection as well.
\hfill $\Box$
\begin{cla}\label{TB041}
If $\la \k _i :i\in I\ra$ is a sequence of cardinals  and some of them is infinite, then
\begin{equation}\label{EQB033}
\la \k _i :i\in I\ra \mbox{ is reversible } \Leftrightarrow \la \k _i :i\in I\ra \mbox{ is finite-to-one}.
\end{equation}
\end{cla}
\dok
Let $i^*\in I$, where $\k_{i^*} \geq \o$.
By Claim \ref{TB040} the implication ``$\Rightarrow$" remains to be checked and we prove its contrapositive.
Suppose that $|I_{\k _0}|\geq \o$, for some cardinal $\k_0$.

If $\k_0 \leq \k_{i^*}$, then we choose different $i_n\in I_{\k _0}\setminus \{ i^*\}$, $n\in \o$,
and define a surjection $f:I\rightarrow I$ by:
$$
f (i)=\left\{  \begin{array}{cl}
                 i^*,          & \mbox{ if } i\in \{ i^*, i_0\} ,\\
                 i_{n-1} ,     & \mbox{ if } i=i_n, \mbox{ and } n\geq 1,  \\
                 i  ,          & \mbox{ if } i\in I \setminus (\{ i^*\} \cup \{ i_n :n\in \o\}).
              \end{array}
     \right.
$$
Now, for $j\in I \setminus (\{ i^*\} \cup \{ i_n :n\in \o\}) $ we have $f^{-1}[\{ j\}]=\{ j\}$;
for $n\in \N$ we have $f^{-1}[\{ i_{n-1}\}]=\{ i_{n}\}$ and $\k _{i_n}=\k _{i_{n-1}}=\k _0$;
finally $f^{-1}[\{ i^*\}]=\{ i^*, i_0\}$ and $\k _{i^*}=\k _{i^*}+ \k _0=\k _{i^*}+ \k _{i_0}$.
So (\ref{EQB032}) is true and, since $f$ is not a bijection,
the sequence $\la \k _i :i\in I\ra$ is not reversible.

If $\k_0 > \k_{i^*}$, then
we choose different $i_n\in I_{\k _0}$, for $n\in \o$,
and define a non-injective surjection $f:I\rightarrow I$ by:
$$
f (i)=\left\{  \begin{array}{cl}
                 i_0,          & \mbox{ if } i\in \{ i_0, i_1\} ,\\
                 i_{n-1} ,     & \mbox{ if } i=i_n, \mbox{ and } n\geq 2,  \\
                 i  ,          & \mbox{ if } i\in I \setminus  \{ i_n :n\in \o\}.
              \end{array}
     \right.
$$
Since $f^{-1}[\{ i_0\}]=\{ i_0, i_1\}$ and $\k_0$ is an infinite cardinal, we have
$\k _{i_0}=\k _0 =\k _{0}+\k _{0}=\k _{i_0}+\k _{i_1}$.
So (\ref{EQB032}) is true and
$\la \k _i :i\in I\ra$ is not reversible again.
\hfill $\Box$
\subsection{Reversible sequences of natural numbers}
Here we characterize reversible sequences of the form $\la n_i :i\in I\ra \in {}^{I}\N$, where $I\neq \emptyset$.
Clearly, $I=\bigcup _{m\in \N}I_m$, where
$$
I_m =\{ i\in I : n_i =m\},\;  \mbox{ for } m\in \N ,
$$
and the following statement is the main result of this paragraph.
\begin{prop}\label{TB037}
A sequence $\la n_i :i\in I\ra \in {}^{I}\N$ is reversible if and only if the set $K:=\{ m\in \N : |I_m|\geq \o\}$ is independent
and, if $K$ is a non-empty set, 
then at most finitely many elements of the set $\{ n_i :i\in I \}$ are divisible by the $\gcd (K)$.
\end{prop}
A proof of Proposition \ref{TB037} is given in the sequel.
First for $d\in \N$ we define $d\N :=\{ dk:k\in \N\}$ and recall some  facts from elementary number theory (giving their proofs for reader's convenience).
\begin{fac}\label{TB032}
Let $K$ be a nonempty subset of $\N$ and $d=\gcd (K)$. Then we have:

(a) If $|K|=\o $, then $\gcd(K^\prime)=d$, for some finite $K^\prime\subset K$;

(b) If $d=1$, then there is $M\in \N$ such that $[M,\infty )\subset \la K \ra$;

(c) If $d>1$,  then there is $M\in \N$ such that $[dM ,\infty )\cap d\N\subset \la K \ra\subset d\N$;

(d) Each independent set is finite.
\end{fac}
\dok
(a) Let $K=\{ n_r:r\in \N\}$ and $d_r =\gcd \{ n_1,\dots ,n_r\}$, for $r\in \N$. Then $d_1 \geq d_2 \geq \dots$
and, hence, there is $s\in \N$ such that $d_r=d_s$, for all $r\geq s$. Clearly we have
$d\leq d_{s}$ and, since $d_{s}$ divides all $n_r$'s,  $d\geq d_{s}$, by the maximality of $d$.
Now we take $K'=\{ n_1,\dots ,n_{s}\}$.

(b) By (a) there is $K^\prime=\{ n_1,\dots ,n_{s}\}\subset K$ such that $\gcd(K^\prime)=1$.
By B\'ezout's lemma 
there are $a_r\in\mathbb{Z}$, for $1\leq r\leq s$, such that
$\sum_{r=1}^s a_r n_r=1$, which for $M:=n_1\sum_{r=1}^s|a_r|n_r$, and for any $m\in\{0,1,\ldots,n_1-1\}$,
implies $M+m=\sum_{r=1}^s(n_1|a_r|+m a_r)n_r\in\langle K^\prime\rangle$; so, $[M,M+n_1)\subset\langle K^\prime\rangle$.
Since $k n_1\in\langle K^\prime\rangle$, we also have that $[M+k n_1,M+(k+1)n_1)\subset\langle K^\prime\rangle$, for any $k\in\mathbb{N}$. Hence, $[M,\infty)\subset\langle K^\prime\rangle\subset\langle K\rangle$.

(c) It is clear that $\la K \ra\subset d\N$.
By (a) there is $K^\prime=\{ n_1,\dots ,n_{s}\}\subset K$ such that $\gcd(K^\prime)=d$
and, hence, $K^\prime=\{d m_1, \dots ,d m_s\}$, where $\gcd(\{m_1 ,\dots ,m_s\})=1$.
By (b) there is $M\in\mathbb{N}$ such that $[M,\infty)\subset\langle\{m_1, \dots ,m_s\}\rangle$, so $[dM,\infty)\cap d\N\subset\langle K^\prime\rangle\subset\langle K\rangle$.

(d) If $K$ is an infinite set, then by (a) there is a finite $K^\prime\subset K$ such that $\gcd(K^\prime)=\gcd(K)=d$. Since $K\setminus K^\prime\subset d\N$ is infinite, for every $M\in\mathbb{N}$  we have $(K\setminus K^\prime)\cap[dM,\infty)\cap d\N\neq\emptyset$. By (c) there is $M\in\mathbb{N}$ such that $[dM,\infty)\cap d\N\subset\langle K^\prime\rangle$. Then $(K\setminus K^\prime)\cap\langle K^\prime\rangle\supset(K\setminus K^\prime)\cap[dM,\infty)\cap d\N\neq\emptyset$. Take $n\in(K\setminus K^\prime)\cap\langle K^\prime\rangle$. Then $n\in K$ and $n\in\langle K^\prime\rangle\subset\langle K\setminus\{n\}\rangle$, which means that the set $K$ is not independant.
\hfill $\Box$
\paragraph{Proof of ``$\Rightarrow$" of Proposition \ref{TB037}}
Let $\la n_i :i\in I\ra$  be a reversible sequence.

First, suppose that the set $K$ is not independent.
Then for some $m\in K$ there are $s>0$, $k_r \in \N$ and different $m_r\in K\setminus \{ m \}$, for $0\leq r<s$, such that
\begin{equation}\label{EQB013}\textstyle
m=\sum _{0\leq r<s}k_r m_r .
\end{equation}
We take countable subsets with 1-1 enumerations
$$
I'_m =\{ j_l :l\in \o\} \subset I_m
$$
$$
I'_{m _r} =\{ i^r_l :l\in \o\} \subset I_{m_r}, \mbox{ for } r<s,
$$
and define $f : I\rightarrow I$ by
$$
f (i)=\left\{  \begin{array}{cl}
                 j_0,          & \mbox{ if } i=i^r_l, \mbox{ where } r<s \mbox{ and } l<k_r ,\\
                 i^r_{l-k_r} , & \mbox{ if } i=i^r_l, \mbox{ where } r<s \mbox{ and } l\geq k_r,  \\
                 j_{l+1},      & \mbox{ if } i=j_l ,  \mbox{ where }  l\in \o , \\
                 i  ,          & \mbox{ if } i\in I \setminus (I'_m \cup \bigcup _{r<s} I'_{m_r}).
              \end{array}
     \right.
$$
It is easy to see that $f [I'_m \cup \bigcup _{r<s} I'_{m_r}] =I'_m \cup \bigcup _{r<s} I'_{m_r}$ so $f$
is a surjection, satisfies (\ref{EQB004}) and it is not 1-1, which gives a contradiction.
So the set $K$ is independent and, by Fact \ref{TB032}(d), $|K |<\o$.

Second, suppose that  $K\neq\emptyset$, $d=\gcd (K)$ and  $|\{ n_i :i\in I \} \cap d\N | =\o$.
\begin{cla}\label{TB046}
There is a sequence $\la q_r :r\in \o\ra$ in $\{ n_i :i\in I\} \cap \la K\ra \setminus K$
such that
\begin{equation}\label{EQB021}
\forall r\in \o \;\; q_{r+1}-q _r \in \la K\ra.
\end{equation}
\end{cla}
\dok
Since $K$ is a finite set, by Fact \ref{TB032}(c)
there is $M\in \N$ such that $M>\max K$ and
\begin{equation}\label{EQB022}
\la K\ra \cap [dM ,\infty )=d\N \cap [dM ,\infty )=\{ dm: m\geq M\}.
\end{equation}
So $\{ n_i :i\in I \} \cap d\N \cap [dM ,\infty )=\{ n_i :i\in I \} \cap \la K\ra \cap [dM ,\infty )$ is an infinite set.
Let $\{ n_i :i\in I \} \cap \la K\ra \cap [dM ,\infty )=\{ n_{i_k}: k\in \o\}$,
where $n_{i_0}< n_{i_1}< n_{i_2}<\dots$.
By recursion we easily construct a sequence $\la k_r :r\in \o\ra$ in $\o$
such that $n_{i_{k_{r+1}}} - n_{i_{k_r}} \geq dM$, which implies that $n_{i_{k_r}}\in \la K \ra\setminus K$ and
$n_{i_{k_{r+1}}} - n_{i_{k_r}} \in \la K \ra$. Defining $q_r=n_{i_{k_r}}$, for $r\in \o$, we finish the proof of Claim \ref{TB046}.
\kdok
For $r\in \o$ we choose $i_r\in I$ such that
\begin{equation}\label{EQB024}
q_r=n_{i_r}\in\la K\ra \setminus K.
\end{equation}
Then by (\ref{EQB021}) and (\ref{EQB024}), $\{ I_m : m\in K \} \cup \{ I_{n_{i_r}}: r\in \o \}$ is a family of pairwise disjoint
subsets of $I$.
For each $m\in K$ we choose a countably infinite, co-infinite subset $I'_m$ of $I_m$ and an 1-1 enumeration of $I'_m$, that is
\begin{equation}\label{EQB025}
I'_m =\{ i^m_l:l\in \o \}\subset I_m \;\;\land \;\;|I'_m|=\o \;\;\land \;\; |I_m \setminus I'_m|\geq \o,
\end{equation}
and in this way we obtain an ``one-to-one matrix indexing" $\{ i^m_l : \la m,l\ra \in K\times \o\}$
of the set $\bigcup _{m\in K}I'_m$.

Now, by (\ref{EQB021}), (\ref{EQB024}) and since the sets $I'_m$ are infinite, we can choose non-empty sets $L_r$, for $r\in \o$, such that

(l1) $L_r \in [K\times \o]^{<\o}$,

(l2) $r_1\neq r_2 \Rightarrow L_{r_1} \cap L_{r_2}=\emptyset$,

(l3) $q_0=n_{i_0}=\sum _{\la m,l\ra\in L_0}n_{i^m_l}$,

(l4) $q_{r+1}-q_r = n_{i_{r+1}}-n_{i_r}=\sum _{\la m,l\ra\in L_{r+1}}n_{i^m_l}$, for $r\in \o$.

\noindent
First,  defining for each $r\in \o$

(g1) $g(i_r)=i_{r+1}$,

(g2) $g(i^m_l)=i_r$, for all $\la m,l \ra\in L_r$,

\noindent
by (l2) we obtain a surjection
\begin{equation}\label{EQB026}\textstyle
g: \{ i^m_l : \la m,l \ra\in \bigcup _{r\in \o }L_r \} \cup \{ i_r :r\in \o\}\rightarrow \{ i_r :r\in \o\}.
\end{equation}
Since $g^{-1}[\{ i_0\}]=\{ i^m_l:\la m,l \ra\in L_0\}$ by (l3) we have
\begin{equation}\label{EQB027}\textstyle
n_{i_0}=\sum _{\la m,l\ra\in L_0}n_{i^m_l}= \sum _{i\in g^{-1}[\{ i_0\}]}n_i .
\end{equation}
Since $g^{-1}[\{ i_{r+1}\}]=\{ i_r\} \cup \{ i^m_l:\la m,l \ra\in L_{r+1}\}$ by (l4) we have
\begin{equation}\label{EQB028}\textstyle
n_{i_{r+1}}=n_{i_r}+\sum _{\la m,l\ra\in L_{r+1}}n_{i^m_l}= \sum _{i\in g^{-1}[\{ i_{r+1}\}]}n_i .
\end{equation}
By (\ref{EQB024}) we have $n_{i_0}\not\in K$ so, by (\ref{EQB027}) we have $|L_0|>1$ and, hence, $g$ is a surjection but not a bijection.
In addition, by (\ref{EQB027}) and (\ref{EQB028})
\begin{equation}\label{EQB029}\textstyle
\forall j\in \{ i_r :r\in \o\} \;\; n_{j}= \sum _{i\in g^{-1}[\{ j\}]}n_i .
\end{equation}
For each $m\in K$ we have $I_m \cap \{ i^{m'}_{l'} : \la m',l' \ra\in \bigcup _{r\in \o }L_r \}\subset I_m '$
so by (\ref{EQB025}) we have $|I_m|=|I_m \setminus \{ i^{m'}_{l'} : \la m',l' \ra\in \bigcup _{r\in \o }L_r \} |$ and, hence,
there are bijections
\begin{equation}\label{EQB030}\textstyle
g_m : I_m \setminus  \{ i^{m'}_{l'} : \la m',l' \ra\in \bigcup _{r\in \o }L_r \} \rightarrow I_m.
\end{equation}
So, for $j\in I_m$ we have $g_m ^{-1}[\{ j\}]=\{ i_j\}$, for some $i_j\in \dom g_m$ and, since $i,i_j\in I_m$,
\begin{equation}\label{EQB031}\textstyle
\forall j\in I_m \;\;n_j =n_{i_j}= \sum _{i\in g_m^{-1}[\{ j\}]}n_i .
\end{equation}
By (\ref{EQB026}) and (\ref{EQB030}) the function $g\cup \bigcup _{m\in K}g_m$ maps the set
$\bigcup _{m\in K}I_m \cup \{ i_r :r\in \o\}$ onto itself and, defining
$$\textstyle
f=g\cup \bigcup _{m\in K}g_m \cup \id _{I\setminus (\bigcup _{m\in K}I_m \cup \{ i_r :r\in \o\})}
$$
by (\ref{EQB029}) and (\ref{EQB031}) we obtain a surjection $f:I\rightarrow I$ which is not a bijection and satisfies (\ref{EQB004}),
which contradicts our assumption that the sequence $\la n_i :i\in I\ra$ is reversible.
The implication ``$\Rightarrow$" of Proposition \ref{TB037} is proved.
\hfill $\Box$
\paragraph{Proof of ``$\Leftarrow$" of Proposition \ref{TB037}}
Let $K$ be an independent set and, if $K\neq\emptyset$, let $|\{ n_i :i\in I \} \cap d\N |<\o$, where $d=\gcd (K)$.

Suppose that the sequence $\la n_i :i\in I\ra$ is not reversible.
Then by Claim \ref{TB040} we have $K\neq\emptyset$ and, hence, $|\{ n_i :i\in I \} \cap d\N |<\o$.
Let $f:I\rightarrow I$ be a surjection such that
\begin{equation}\label{EQB004}\textstyle
\forall j\in I \;\; n_j=\sum _{i\in f^{-1}[\{ j \}]}n_i.
\end{equation}
\begin{equation}\label{EQB005}\textstyle
J:= \{ j\in I : |f^{-1}[\{ j \}]|>1\} \neq \emptyset.
\end{equation}

\begin{cla}\label{TB031}
(a) For each $i\in I$ we have $n_i\leq n_{f(i)}$.

(b) For each $j\in I$ there is a sequence
$\la i^j _k :k\in \N\ra$ in $I$ such that
\begin{equation}\label{EQB007}\textstyle
f(i^j_1)=j \;\;\land \;\; \forall k\in \N \;\; f(i^j_{k+1})=i^j_k ,
\end{equation}
\begin{equation}\label{EQB008}\textstyle
\dots n_{i^j_{k+1}}\leq n_{i^j_k} \leq \dots n_{i^j_{3}}\leq n_{i^j_2} \leq n_{i^j_{1}}\leq n_{j}.
\end{equation}

(c) If, in addition, $n_{i^j_{1}}< n_{j}$ in (\ref{EQB008}), then $i^j_k\neq i^j_l$, whenever $k\neq l$.
\end{cla}
\dok
(a) follows from (\ref{EQB004}).

(b) If $j\in I$,
then, since $f$ is an onto mapping, there is $i^j_1\in I$ such that $f(i^j_1)=j$,
there is $i^j_2\in I$ such that $f(i^j_2)=i^j_1$,
there is $i^j_3\in I$ such that $f(i^j_3)=i^j_2$,
and so on.
So in this way we obtain a sequence $\la i^j _k :k\in \N\ra \in {}^{\N}I$
satisfying (\ref{EQB007})
which, together with (a), gives (\ref{EQB008}).

(c) If $n_{i^j_{1}}< n_{j}$ then, by (\ref{EQB008}), $n_{i^j_{k}}< n_{j}$, for all $k\in \N$ and, hence,
\begin{equation}\label{EQB010}\textstyle
\forall k\in \N \;\; i^j_k\neq j .
\end{equation}
On the contrary, let $k$ be the minimal element of $\N$ such that $i^j_k= i^j_l$, for some $l>k$.
Then by (\ref{EQB007}), for $k=1$ we would have  $i^j_{l-1}= f(i^j_l)=f(i^j_k)=f(i^j_1)=j$, which is impossible by (\ref{EQB010}).
For $k>1$ we would have $i^j_{l-1}= f(i^j_l)=f(i^j_k)=i^j_{k-1}$, which is false by the minimality of $k$.
\hfill $\Box$
\begin{cla}\label{TB029}
There is a sequence $\la p_r :r\in \o \ra$ in $\N$
such that, defining for convenience $p_{-1}:=0$, for each $r\in \o$ we have:

(i) $p_r=\min \{ n_j : j\in J \land n_j >p_{r-1}\}$,

(ii) $\forall j\in I_{p_r}\cap J\;\; \forall i\in f^{-1}[\{ j \}] \;\; n_i \in K \cup \{ p_s : 0\leq s<r\}$,

(iii) $p_r\in \la K\ra\setminus K$,

(iv) $\exists i\in I_{p_r} \;\; ( f(i)\in J \land n_{f(i)}>p_r)$,

(v) $\{ n_j : j\in J\}\cap [1,p_r]=\{ p_s : 0\leq s\leq r\}$.
\end{cla}
\dok
We construct the sequence by recursion.

First, by (\ref{EQB005}) we have $J\neq \emptyset$ so
$\emptyset \neq \{ n_j : j\in J\}=\{ n_j : j\in J\land n_j >0\}\subset \N$ and defining
\begin{equation}\label{EQB011}\textstyle
p_0=\min \{ n_j : j\in J\}
\end{equation}
we see that the sequence $\la p_0\ra$ satisfies (i).

(ii) Let $j\in I_{p_0}\cap J$ and $i\in f^{-1}[\{ j \}]$.
Then, since $j\in J$, by (\ref{EQB005}) we have $|f^{-1}[\{ j \}]|>1$
and, by (\ref{EQB004}), $n_j =\sum _{i'\in f^{-1}[\{ j \}]}n_{i'}$, so $n_i <n_{j}$.
As in Claim \ref{TB031} we define $i^{j}_k\in I$, for $k\in \N$, satisfying $i^{j}_1:=i$, (\ref{EQB007}) and (\ref{EQB008})
and so we obtain $\dots n_{i^{j}_{3}}\leq n_{i^{j}_2} \leq n_{i^{j}_{1}}< n_{j} $.
Assuming that $n_{i^{j}_{k+1}}< n_{i^{j}_k}$ for some $k\in \N$, since $f(i^{j}_{k+1})=i^{j}_k$ by (\ref{EQB004})
we would  have $i^{j}_k\in J$ and $n_{i^{j}_k}<n_{j}=p_0$, which is, by (\ref{EQB011}), impossible.
Thus there is $m\in \N$ such that $n_{i^{j}_k}=m$, for all $k\in \N$.
By Claim \ref{TB031}(c) we have $i^{j}_k\neq i^{j}_l$, whenever $k\neq l$, thus $|I_m|\geq \o$. So $n_i =n_{i^{j}_1}=m\in K$.

(iii) By the previous item and (\ref{EQB004}) we have $p_0=n_{j} \in \la K\setminus \{ p_0\}\ra$ and, since the set  $K$ is independent, $p_0\not \in K$.

(iv) By (iii) we have $p_0\not\in K$, that is $|I_{p_0}|<\o$. Suppose that $f[I_{p_0}]\subset I_{p_0}$.
Then by (\ref{EQB004}) $f\upharpoonright I_{p_0}$ is an injection and, since the set $I_{p_0}$ is finite, $f[I_{p_0}]= I_{p_0}$.
By (\ref{EQB011}) there is $j\in I_{p_0}\cap J$ and by the previous conclusion, $j=f(i)$, for some $i\in I_{p_0}$, which implies that
$n_i=n_j=p_0$. But this contradicts the fact that $j\in J$. So, there is $i\in I_{p_0} $ such that $f(i)\not\in I_{p_0}$
and, hence, $n_{f(i)}> n_i=p_0$ and $f(i)\in J$.

(v) By (\ref{EQB011}) we have $\{ n_j : j\in J\}\cap [1,p_0]=\{ p_0\}$.

Suppose that $\la p_0, \dots ,p_r\ra$ is a sequence satisfying (i)--(v). By (iv) there is $j\in J$ such that $n_j >p_r$
and defining
\begin{equation}\label{EQB012}\textstyle
p_{r+1}=\min \{ n_j : j\in J \land n_j >p_r\} .
\end{equation}
we have (i).

(ii) Let $j\in I_{p_{r+1}}\cap J$ and $i\in f^{-1}[\{ j \}]$. Then, since $j\in J$, by (\ref{EQB005}) we have $|f^{-1}[\{ j \}]|>1$
and, by (\ref{EQB004}), $n_j =\sum _{i'\in f^{-1}[\{ j \}]}n_{i'}$, so $n_i <n_{j}$.
Again, as in Claim \ref{TB031} we define $i^{j}_k\in I$, for $k\in \N$, satisfying $i^{j}_1:=i$, (\ref{EQB007}) and (\ref{EQB008})
and so we obtain $\dots n_{i^{j}_{3}}\leq n_{i^{j}_2} \leq n_{i^{j}_{1}}< n_{j} $.

If $n_{i^{j}_{k+1}}< n_{i^{j}_k}$ for some $k\in \N$, let $k$ be the minimal such $k$.
Then
\begin{equation}\label{EQB014}
n_{i^{j}_{k+1}}<n_{i^{j}_k} = \dots = n_{i^{j}_2} = n_{i^{j}_{1}}=n_i< n_{j}=p_{r+1} .
\end{equation}
In addition, since $f(i^{j}_{k+1})=i^{j}_k$, by (\ref{EQB004}) we have $i^{j}_k\in J$
which implies that $n_{i^{j}_k}\in \{ n_j:j\in J\} \cap [1,p_{r+1})$ and, by (\ref{EQB012}),
$n_{i^{j}_k}\in \{ n_j:j\in J\} \cap [1,p_r]$.
So, by (v), there is $s_0\leq r$ such that $n_{i^{j}_k}=p_{s_0}$
and, by (\ref{EQB014}), $n_i=p_{s_0}\in \{ p_s : 0\leq s<r+1\}$.

Otherwise, there is $m\in \N$ such that $n_{i^{j}_k}=m$, for all $k\in \N$.
By Claim \ref{TB031}(c) we have $i^{j}_k\neq i^{j}_l$, whenever $k\neq l$,
thus $|I_m|\geq \o$,
and, hence, $m\in K$.
So $n_i =n_{i^{j}_1}=m\in K$ and (ii) is true indeed.

(iii) By (\ref{EQB012}) there is $j\in J$ such that $p_{r+1}=n_j>p_r$.
Thus $j\in I_{p_{r+1}}\cap J$ and, by (ii) and (\ref{EQB004}), $n_j$
is a sum of at least two integers from $K\cup \{ p_s : 0\leq s\leq r\}$.
By (iii) of the induction hypothesis we have $p_s\in \la K\ra$, for $0\leq s\leq r$, and, hence, $p_{r+1}\in \la K\setminus \{p_{r+1} \}\ra$.
Since the set $K$ is independent we have $p_{r+1}\not\in  K$.

(iv) Since $p_{r+1}\not\in K$ we have $|I_{p_{r+1}}|<\o$. Suppose that $f[I_{p_{r+1}}]\subset I_{p_{r+1}}$.
Then by (\ref{EQB004}) $f\upharpoonright I_{p_{r+1}}$ is an injection and, since the set $I_{p_{r+1}}$ is finite, $f[I_{p_{r+1}}]= I_{p_{r+1}}$.
By (\ref{EQB012}) there is $j\in I_{p_{r+1}}\cap J$ and, since $f[I_{p_{r+1}}]= I_{p_{r+1}}$, $j=f(i)$, for some $i\in I_{p_{r+1}}$, which implies that
$n_i=n_j=p_{r+1}$. But this contradicts the fact that $j\in J$. So, there is $i\in I_{p_{r+1}} $ such that $f(i)\not\in I_{p_{r+1}}$
and, hence, $n_{f(i)}> n_i=p_{r+1}$ and $f(i)\in J$.

(v) By (\ref{EQB012}) and the induction hypothesis we have $\{ n_j : j\in J\}\cap [1,p_{r+1}]=\{ p_s : 0\leq s\leq r+1\}$.
Thus the recursion works.
\kdok
Now, by Claim \ref{TB029}(v), (iii) and (i),
$\{ n_j :j\in J\} =\{ p_r :r\in \o\} \subset \la K\ra \setminus K$ and $p_0 <p_1 <\dots < p_r < \dots $,
which implies that $| \{ n_i :i\in I \}\cap\la K \ra |=\o$. Since, by Fact \ref{TB032}(c), $\la K\ra \subset d\N$,
we have $|\{ n_i :i\in I \} \cap d\N | =\o$ and we obtain a contradiction.
\hfill $\Box$
\paragraph{Reversible functions in the Baire space}
Each countable sequence of natural numbers $\la n_i :i\in \N\ra \in {}^{\N}\N$ can be regarded
as a function $\f :\N \rightarrow \N$, where $\f(i)=n_i$, for $i\in \N$, and, hence, as an element
of the Baire space $\N^\N$ with the standard topology (see \cite{Kech}). So we can consider the set of reversible functions
belonging to $\N^\N$,
$$
\Nrev \!:= \Big\{ \f \in \N^\N : \neg \exists f\in \Sur(\N)\setminus\Sym(\N ) \;\forall j\in \N \;\;\f (j)=\!\sum _{i\in f^{-1}[\{ j \}]} \f (i)  \Big\} .
$$
\begin{te}
$\Nrev $ is a dense $F_{\sigma\delta\sigma}(=\Sigma ^0_4)$ subset of $\N^\N$ of size ${\mathfrak c}$.
\end{te}
\dok
If $B=\bigcap _{k\leq n}\pi ^{-1}_{i_k}[\{ j_k\}]$ is a basic open set, then, since the finite function $p=\{ \la i_k ,j_k\ra : k\leq n\}$
can be extended to an finite-to-one function $\f \in \N^\N$ and  by Proposition \ref{TB037} we have $\f\in \Nrev$, it follows that
$B\cap \Nrev \neq \emptyset$ so $\Nrev $ is  dense in $\N^\N$.  $|\Nrev|={\mathfrak c}$ follows from the fact that $\N^\N$ contains
${\mathfrak c}$-many injections.

Let $\I$ be the set of non-empty independent subsets of $\N$ and, for $K\in \I$, let $d_K:=\gcd (K)$.
Then by Proposition \ref{TB037}
\begin{equation}\label{EQB023}\textstyle
\Nrev = A\cup \bigcup _{K\in \I}\;B_K \cap C_K \cap D_K,
\end{equation}
where
\begin{eqnarray*}
A   & := & \Big\{ \f\in \N^\N : \forall m\in \N \;\; ( \f(i)=m \mbox{ for } <\o\mbox{-many } i\in \N)\Big\} , \\
    & =  &\textstyle \bigcap _{m\in \N} \bigcup _{k\in \N} \bigcap _{i\geq k} \pi ^{-1}_i[\N \setminus \{ m\}] , \\
B_K & := & \Big\{ \f\in \N^\N : \forall m\in K \;\; ( \f(i)=m \mbox{ for } \o\mbox{-many } i\in \N)\Big\} ,\\
    & =  &\textstyle \bigcap _{m\in K} \bigcap _{k\in \N} \bigcup _{i\geq k} \pi ^{-1}_i[\{ m\}] , \\
C_K & := & \Big\{ \f\in \N^\N : \forall m\in \N \setminus K \;\; ( \f(i)=m \mbox{ for } <\o\mbox{-many } i\in \N)\Big\} , \\
    & =  &\textstyle \bigcap _{m\in \N\setminus K} \bigcup _{k\in \N} \bigcap _{i\geq k} \pi ^{-1}_i[\N \setminus \{ m\}] , \\
D_K & := & \Big\{ \f\in \N^\N : \f(i)\in d\N \mbox{ for } <\o\mbox{-many } i\in \N\Big\} \\
    & =  &\textstyle \bigcup _{m\in \N} \bigcap _{i\geq m} \bigcap _{k\in \N} \pi ^{-1}_i[\N \setminus \{ dk\}] .
\end{eqnarray*}
So, for $K\in \I$ we have $B_K \in G_\delta$, $D_K\in F_\sigma$ and $C_K\in F_{\sigma\delta}$, which implies that
$B_K \cap C_K \cap D_K \in F_{\sigma\delta}$ and, since by Fact \ref{TB032}(d) we have $\I \subset [\N ]^{<\o}$, it follows that
$\bigcup _{K\in \I}\;B_K \cap C_K \cap D_K \in F_{\sigma\delta\sigma}$. Since $A\in F_{\sigma\delta}\subset F_{\sigma\delta\sigma}$,
by (\ref{EQB023}) we have $\Nrev \in F_{\sigma\delta\sigma}=\Sigma ^0_4$.
\hfill $\Box$
\begin{rem}\label{RB000}\rm
Let the equivalence relation $\sim$ on $\N^\N$ be defined by $\f \sim \p$ iff there is $f\in \Sym (\N)$ such that
$\f=\p\circ f$. It is evident that the set $\Nrev$ is $\sim$-invariant, that is $\p\sim \f\in \Nrev$ implies $\p\in \Nrev$.

But $\Nrev$ is not a subsemigroup of $\la \N ^\N ,\circ \ra$ (it is not closed under composition). Let $\N \setminus \{ 2\}=A \cup B$  and $\N =C\cup D \cup E$ be partitions,
where $A,B,C,D,E \in [\N ]^\o$ and $|A\cap (2\N +1)|=|B\cap (2\N +1)|=\o$.  Then, by Proposition \ref{TB037},
$
\f =\{ \la 2,2 \ra\}\cup (A\times \{ 3 \}) \cup (B\times \{ 5 \})\in \Nrev.
$

If $\p_{DA}:D\rightarrow A\cap (2\N +1)$ and $\p_{EB}:E\rightarrow B\cap (2\N +1)$ are bijections then, by Proposition \ref{TB037} again,
$
\p = (C\times \{ 2\}) \cup \p_{DA}\cup\p_{EB}\in \Nrev .
$
But $\f \circ \p \not\in \Nrev$, because the set $\{ 2,3,5\}$ is not independent.
\end{rem}

\end{document}